\documentclass[11pt]{amsart}
\usepackage[T1]{fontenc}
\usepackage{lmodern,amsmath,amssymb,mathtools}
\usepackage[margin=1.1in]{geometry}
\usepackage[colorlinks=true,linkcolor=blue,citecolor=blue,urlcolor=blue]{hyperref}
\newtheorem{theorem}{Theorem}[section]
\newtheorem{lemma}{Lemma}[section]
\newtheorem{proposition}{Proposition}[section]
\newtheorem{corollary}{Corollary}[section]
\theoremstyle{remark}
\newtheorem{remark}{Remark}[section]
\newcommand{\Ric}{\operatorname{Ric}}
\newcommand{\Vol}{\operatorname{Vol}}
\newcommand{\Crit}{\operatorname{Crit}}
\newcommand{\dd}{\,d v_g}
\newcommand{\ds}{\,d\sigma_g}

\title{Scalar curvature bounds and rigidity of Ricci shrinkers}
\author{Jianyu Ou}
\address{School of Mathematical Sciences, Xiamen University, Xiamen 361005, China}
\email{oujianyu@xmu.edu.cn}
\author{Guoqiang Wu}
\address{School of Science, Zhejiang Sci-Tech University, Hangzhou 310018, China}
\email{gqwu@zstu.edu.cn}
\date{}
\begin{document}
\begin{abstract}
We prove that a complete noncompact gradient shrinking Ricci soliton of dimension at least two is a finite quotient of an Einstein cylinder if its scalar curvature is bounded below by the cylindrical value. The corresponding pointwise strict bound implies compactness. We also establish one-endedness criteria using curvature information at critical points of the potential, including degenerate ones, and obtain sharp one-endedness and rigidity results under global Ricci norm bounds.
\end{abstract}
\maketitle

\section{Introduction}\label{sec:intro}
Curvature bounds on a complete gradient shrinking Ricci soliton provide information about its compactness and geometry at infinity. Two natural questions are whether the soliton is compact and, in the noncompact case, how many ends it has. In this paper, we study these questions under scalar curvature bounds. We also establish criteria for one-endedness using Ricci curvature information at the critical points of the potential, and prove a sharp cylindrical alternative under a global Ricci norm bound.

A complete Riemannian manifold $(M^n,g)$ is called a gradient shrinking Ricci soliton if there exists a smooth function $f$ on $M$ such that, after scaling the metric,
\[
\Ric+\nabla^2f=\frac12g.
\]
Here $\Ric$ denotes the Ricci tensor, $\nabla^2f$ is the Hessian of $f$, and $f$ is called a potential function. Gradient Ricci solitons generalize Einstein metrics, generate self-similar solutions to the Ricci flow, and arise as models for singularity formation. We refer to \cite{Cao} for a nice survey on the subject.

Throughout the paper, we assume that the soliton manifold $M$ is connected and orientable, with dimension $n\ge2$, and normalize the potential by
\[
R+|\nabla f|^2=f,
\]
where $R$ denotes the scalar curvature. All manifolds are smooth and without boundary. Tensor norms are the pointwise Hilbert--Schmidt norms induced by $g$. The one-dimensional complete shrinking soliton is the Gaussian line, with $R=|\Ric|=0$ and two ends; we exclude this elementary case from the statements below.

Several compactness results are known for shrinking Ricci solitons. Munteanu and Wang \cite{MWcompact} proved compactness under nonnegative sectional curvature and positive Ricci curvature. Under the scalar curvature lower bound $R\ge n/2$, a complete gradient shrinking Ricci soliton is compact and Einstein, as proved by Pigola, Rimoldi and Setti and by Mai and Ou \cite{PRS,MO}. It is therefore natural to ask how far the scalar curvature lower bound can be weakened while still implying compactness.

The volume estimates in \cite{CZvolume,Zhang,MWanalysis} imply $\inf_M R\le(n-1)/2$ for every complete noncompact shrinker. This threshold is attained by Einstein cylinders. Fern\'andez-L\'opez and Garc\'ia-R\'io \cite{FLGR} proved rigidity when $R\equiv(n-1)/2$.

Sun and Zhu \cite{SZminimal} obtained splitting under an additional area-minimizing hypersurface assumption, while Chai and Sun \cite{CSspectral} developed warped $\mu$-bubble methods under spectral curvature bounds. Our first main result shows that $R\ge(n-1)/2$ alone forces equality in the noncompact case.

\begin{theorem}\label{thm:main}
Let $(M^n,g,f)$, $n\ge2$, be a connected, complete, noncompact gradient shrinking Ricci soliton satisfying
\[
\Ric+\nabla^2f=\frac12g,\qquad R+|\nabla f|^2=f,
\]
where $R$ is the scalar curvature of $g$.
If
\[
R\ge\frac{n-1}{2}\qquad\text{on }M,
\]
then
\[
R\equiv\frac{n-1}{2}.
\]
Consequently, $(M,g)$ is isometric to a finite quotient of $(N^{n-1}\times\mathbb R,g_N+ds^2)$, where $N$ is a compact Einstein manifold with $\Ric_{g_N}=\tfrac12g_N$.
\end{theorem}

The geometric conclusion follows from the constant scalar curvature rigidity theorem in \cite{FLGR}. If $M$ is simply connected, it is a global product. The Einstein cylinders show that the lower bound in Theorem~\ref{thm:main} is sharp for $n\ge3$.

As a consequence, we obtain compactness when the scalar curvature is strictly above the cylindrical threshold at every point. Unlike the direct consequence of the volume estimates, this criterion requires no uniform positive gap.

\begin{corollary}\label{cor:compact}
Let $(M^n,g,f)$ be a connected, complete gradient shrinking Ricci soliton satisfying
\[
\Ric+\nabla^2f=\frac12g,\qquad R+|\nabla f|^2=f,
\]
where $R$ is the scalar curvature of $g$.
If
\[
R(x)>\frac{n-1}{2}\qquad\text{for every }x\in M,
\]
then $M$ is compact.
\end{corollary}

Scalar curvature bounds also give information about the ends of a shrinker. Munteanu, Schulze and Wang \cite{MSW} proved finiteness of the number of ends under an integral condition on the scalar curvature; in particular, their result applies when the scalar curvature has quadratic decay. In the K\"ahler case, Munteanu and Wang \cite{MWkahler} proved that every complete noncompact shrinking soliton has one end. For general shrinking Ricci solitons of dimension $n\ge4$, they proved that $R\le n/3$ implies that $M$ has one end \cite{MWends}.

Combining the scalar upper bound criterion for one-endedness with Theorem~\ref{thm:main}, we obtain the following restriction on the scalar curvature of a shrinker with more than two ends.

\begin{corollary}\label{cor:range}
Let $(M^n,g,f)$, $n\ge4$, be a connected, complete, noncompact gradient shrinking Ricci soliton satisfying
\[
\Ric+\nabla^2f=\frac12g,\qquad R+|\nabla f|^2=f,
\]
where $R$ is the scalar curvature of $g$.
If $M$ has more than two ends, then
\[
\sup_MR>\frac n3,\qquad\inf_MR<\frac{n-1}{2}.
\]
\end{corollary}
Indeed, the first inequality follows from Theorem~\ref{thm:MWends}. If the second inequality failed, Theorem~\ref{thm:main} would imply that $M$ is a finite quotient of an Einstein cylinder, which has at most two ends.

Qu and Wu \cite{QW} established one-endedness under the bound
\[
|\Ric|<\frac{n-2}{2\sqrt n},\qquad n\ge3.
\]
The critical-point argument below improves this strict global threshold to $\sqrt{n-1}/2$. For the non-strict global bound, the only alternative to one-endedness is an Einstein cylinder.

Bertellotti and Buzano \cite{BBends,BBstructure} studied bounds on the number of ends and the geometry at infinity of Ricci shrinkers. He \cite{He} obtained end dichotomies under an Einstein-tensor bound using weighted cohomology.

We next study one-endedness through the Hessian of the potential at its critical points. The following theorem uses curvature information only on the critical set, independently of the scalar lower bound in Theorem~\ref{thm:main}.

\begin{theorem}\label{thm:criticalintro}
Let $(M^n,g,f)$, $n\ge2$, be a connected, complete, noncompact gradient shrinking Ricci soliton satisfying
\[
\Ric+\nabla^2f=\frac12g,\qquad R+|\nabla f|^2=f,
\]
where $R$ is the scalar curvature of $g$.
Let $\Crit(f)=\{x\in M:\nabla f(x)=0\}$ be the critical set of $f$. Suppose that, at every $p\in\Crit(f)$,
\[
\Ric_p\ge0,\qquad R(p)<\frac{n-1}{2}.
\]
Then $M$ has exactly one end, and $H_{n-1}(M;\mathbb Z_2)=0$.
\end{theorem}
Both curvature conditions are imposed only on the critical set. In particular, the Ricci tensor may have negative eigenvalues elsewhere, and the potential need not be Morse or Morse--Bott.

A global Ricci norm bound yields a related dichotomy. At the critical threshold, a shrinker either has one end or is a global Einstein cylinder with exactly two ends.

\begin{theorem}\label{thm:normintro}
Let $(M^n,g,f)$, $n\ge2$, be a connected, complete, noncompact gradient shrinking Ricci soliton satisfying
\[
\Ric+\nabla^2f=\frac12g,\qquad R+|\nabla f|^2=f,
\]
where $R$ is the scalar curvature of $g$.
If
\[
|\Ric|^2\le\frac{n-1}{4}\qquad\text{on }M,
\]
then either $M$ has exactly one end, or $n\ge3$ and $(M,g)$ is a global product $(N^{n-1}\times\mathbb R,g_N+ds^2)$ with exactly two ends, where $N$ is a connected compact Einstein manifold satisfying $\Ric_{g_N}=\tfrac12g_N$. If the norm inequality is strict at one point, then $M$ has exactly one end.
\end{theorem}
In particular, the pointwise strict bound $|\Ric|<\sqrt{n-1}/2$ implies one-endedness without a uniform gap. Product cylinders attain equality, so they must be included under the non-strict bound. A finite cylinder quotient may have one end if its deck action exchanges the two ends; the global product conclusion applies when $M$ has more than one end.

To prove Theorem~\ref{thm:criticalintro}, we first use $f=R$ on $\Crit(f)$ to show that the critical set is compact. The curvature hypotheses imply that the Hessian has two strictly positive eigenvalues at each critical point. A relative Morse perturbation preserves these positive directions, so the resulting handle decomposition has no handles of index $n-1$ or $n$. Mod-two duality then gives a connected boundary for a compact core, and hence one-endedness. For Theorem~\ref{thm:normintro}, a maximum principle argument first yields $R\le(n-1)/2$. The equality case gives rigidity, while the strict case reduces to the same critical-point argument.

The same threshold also gives rigidity when imposed as a lower bound for the Ricci norm. The next theorem characterizes the noncompact case and yields compactness if the lower bound is strict at one point.

\begin{theorem}\label{thm:normlowerintro}
Let $(M^n,g,f)$, $n\ge2$, be a connected, complete gradient shrinking Ricci soliton satisfying
\[
\Ric+\nabla^2f=\frac12g,\qquad R+|\nabla f|^2=f,
\]
where $R$ is the scalar curvature of $g$.
Suppose
\[
|\Ric|^2\ge\frac{n-1}{4}\qquad\text{on }M.
\]
Then either $M$ is compact, or $n\ge3$ and $M$ is a finite quotient of an Einstein cylinder $N^{n-1}\times\mathbb R$, where $N$ is compact and $\Ric_{g_N}=\tfrac12g_N$. In the noncompact case,
\[
R\equiv\frac{n-1}{2},\qquad |\Ric|^2\equiv\frac{n-1}{4}.
\]
If the assumed lower bound is strict at even one point, then $M$ is compact.
\end{theorem}
Theorem~\ref{thm:normlowerintro} follows from the comparison $|\Ric|^2\ge a\ge0\Rightarrow R\ge2a$, obtained by minimizing $R+\varepsilon f$ and letting $\varepsilon\downarrow0$.

Under the complementary bound $R\le(n-1)/2$, we also obtain finite topological type and critical-point criteria for one-endedness, allowing some negative Ricci curvature.

The proof of Theorem~\ref{thm:main} combines a weighted integrability estimate for the scalar curvature excess with a logarithmic barrier. Unless $R\equiv(n-1)/2$, the barrier and the volume lower bound contradict this integrability estimate.

The paper is organized as follows. Section~\ref{sec:prelim} recalls the soliton identities, potential and volume estimates, and regularity and rigidity results used below. Section~\ref{sec:critical} establishes the weighted integrability estimate and the differential inequality for the scalar curvature excess. In Section~\ref{sec:barrier}, we construct the logarithmic barrier and prove Theorem~\ref{thm:main}. Section~\ref{sec:criticalends} proves the critical-set scalar criterion for one-endedness by Morse approximation on a compact core, and Section~\ref{sec:norm} proves the cylindrical alternative under a Ricci norm upper bound. Subsection~\ref{subsec:lowernorm} treats the lower bound, its compactness consequences and sharp infimum estimates. Finally, Subsections~\ref{subsec:scalarupper} and~\ref{subsec:mixedcritical} establish quantitative critical-point obstructions and criteria for one-endedness under the complementary scalar upper bound.

{\bf Acknowledgement:}  The first author is partially supported by NSFC grant No.12401074, 12371061, 12371081.

{\bf Disclosure on AI assistance.} The research problem and the principal mathematical ideas were formulated by the authors. The authors used Astra to assist with checking some computations and polishing the exposition. The authors verified and completed all mathematical arguments and take full responsibility for the content.

\section{Preliminaries}\label{sec:prelim}
We recall some basic facts about gradient shrinking Ricci solitons used in the proof. Unless otherwise stated, $(M^n,g,f)$ is a connected, complete gradient shrinking Ricci soliton. We use the normalization
\begin{equation}\label{eq:normalization}
\Ric+\nabla^2f=\frac12g,\qquad R+|\nabla f|^2=f.
\end{equation}
The second identity is obtained by adding a constant to the potential function. We write
\[
\Delta_f=\Delta-\langle\nabla f,\nabla(\cdot)\rangle
\]
for the weighted Laplacian, and use the notation
\begin{equation}\label{eq:notation}
D_t=\{f<t\},\qquad\delta=\inf_MR.
\end{equation}
The measures $dv_g$ and $d\sigma_g$ denote the Riemannian volume and hypersurface measures, respectively. Constants denoted by $C$ may change from line to line.

\subsection{Basic identities and scalar curvature}

We will use the following soliton identities in the maximum principle and sublevel set arguments; see, for example, \cite{CZpotential,PW}.

\begin{lemma}\label{lem:identities}
Let $(M^n,g,f)$ be a connected, complete gradient shrinking Ricci soliton satisfying
\[
\Ric+\nabla^2f=\frac12g,\qquad R+|\nabla f|^2=f,
\]
where $R$ is the scalar curvature of $g$.
Define $\Delta_f=\Delta-\langle\nabla f,\nabla(\cdot)\rangle$. Then

\begin{equation}\label{eq:fidentities}
\Delta f=\frac n2-R,\qquad\Delta_ff=\frac n2-f,
\end{equation}
and
\begin{equation}\label{eq:Ridentities}
\nabla R=2\Ric(\nabla f,\cdot),\qquad\Delta_fR=R-2|\Ric|^2.
\end{equation}
Here the gradient of $R$ is identified with its metric dual one-form. In particular, if $R$ is constant, then
\[
\Ric(\nabla f,\cdot)=0,\qquad |\Ric|^2=\frac R2.
\]
\end{lemma}

We also record the Gaussian equality case and the compact Einstein rigidity under $R\ge n/2$, the latter established in both \cite{PRS,MO}, for the endpoint arguments.

\begin{lemma}\label{lem:scalar}
Let $(M^n,g,f)$ be a connected, complete gradient shrinking Ricci soliton satisfying
\[
\Ric+\nabla^2f=\frac12g,\qquad R+|\nabla f|^2=f,
\]
where $R$ is the scalar curvature of $g$.
Set $\delta=\inf_M R$. Then $R\ge0$ on $M$, and
\[
0\le\delta\le\frac n2.
\]
If $R$ vanishes at a point, then the soliton is the Gaussian shrinking soliton on $\mathbb R^n$. If $\delta=n/2$, then $M$ is compact, $f$ is constant, and $\Ric=\tfrac12g$.
\end{lemma}

\subsection{The potential function and volume growth}

The quadratic growth estimate of Cao and Zhou \cite[Theorem 1.1]{CZpotential} shows that the potential is a proper exhaustion of a noncompact shrinker. It also compares potential sublevel sets with geodesic balls, allowing us to pass between their volume estimates.

\begin{lemma}\label{lem:potential}
Let $(M^n,g,f)$ be a connected, complete, noncompact gradient shrinking Ricci soliton satisfying
\[
\Ric+\nabla^2f=\frac12g,\qquad R+|\nabla f|^2=f,
\]
where $R$ is the scalar curvature of $g$.
Fix $p\in M$, set $D_t=\{x\in M:f(x)<t\}$, and let $d$ and $B_p(r)$ denote the Riemannian distance and the geodesic ball of radius $r$ centered at $p$, respectively. There exist constants $C,r_0>0$ such that, whenever $d(p,x)\ge r_0$,
\begin{equation}\label{eq:potential}
\frac14(d(p,x)-C)^2\le f(x)\le\frac14(d(p,x)+C)^2.
\end{equation}
In particular, $f$ is proper, attains its minimum, and is unbounded above. After increasing $C$, for all sufficiently large $t$ one has
\begin{equation}\label{eq:inclusions}
B_p(2\sqrt t-C)\subset D_t\subset B_p(2\sqrt t+C).
\end{equation}
\end{lemma}

The following volume estimates will bound the scalar curvature infimum and determine the volume growth in the critical case. The upper bound is due to Zhang \cite{Zhang} when $\delta>0$ and to Cheng--Zhou \cite[Theorem 1.5]{CZvolume}; the lower bound is due to Munteanu--Wang \cite[Theorem 6.1]{MWanalysis}.

\begin{lemma}\label{lem:volume}
Let $(M^n,g,f)$ be a connected, complete, noncompact gradient shrinking Ricci soliton satisfying
\[
\Ric+\nabla^2f=\frac12g,\qquad R+|\nabla f|^2=f,
\]
where $R$ is the scalar curvature of $g$.
Set $\delta=\inf_MR$ and fix $p\in M$. Let $B_p(r)$ denote the geodesic ball of radius $r$ centered at $p$. There are constants $C,r_0>0$ such that, for $r\ge r_0$,
\begin{equation}\label{eq:volumeupper}
\Vol(B_p(r))\le Cr^{n-2\delta},
\end{equation}
and
\begin{equation}\label{eq:volumelower}
\Vol(B_p(r))\ge C^{-1}r.
\end{equation}
Consequently, every complete noncompact gradient shrinking Ricci soliton satisfies $\delta\le(n-1)/2$.
\end{lemma}
The last assertion follows by comparing the exponents of $r$ in \eqref{eq:volumeupper} and \eqref{eq:volumelower}.

\subsection{Analyticity and integration over sublevel sets}

Kotschwar \cite[Corollary 1.3 and the subsequent remark]{Kotschwar} proved that a Ricci soliton metric and its soliton vector field are real analytic in suitable local coordinates. We will use the identity theorem for analytic functions and the fact that the zero set of a nonzero real-analytic function on a connected manifold has measure zero. For a noncompact gradient shrinker, the latter applies to the critical set of the potential.

\begin{lemma}\label{lem:analytic}
Let $(M^n,g)$ be a connected smooth Riemannian manifold without boundary, and let $V$ be a smooth vector field satisfying the Ricci soliton equation
\[
\Ric+\frac12\mathcal L_Vg=\lambda g
\]
for a constant $\lambda\in\mathbb R$, where $\mathcal L_Vg$ denotes the Lie derivative of $g$ along $V$. Then $g$ and $V$ are real analytic in suitable local coordinates. If $V=\nabla f$ for a smooth function $f:M\to\mathbb R$, then $f$ and the scalar curvature $R$ are real analytic in the associated analytic structure. If, in addition, $(M^n,g,f)$ is a complete, noncompact gradient shrinking Ricci soliton satisfying
\[
\Ric+\nabla^2f=\frac12g,\qquad R+|\nabla f|^2=f,
\]
then its critical set
\[
\Crit(f)=\{x\in M:|\nabla f|(x)=0\}
\]
has Riemannian measure zero. Moreover, an analytic function that vanishes on a nonempty open subset of $M$ vanishes identically.
\end{lemma}
Indeed, $df=g(\nabla f,\cdot)$ is analytic, so $f$ is locally analytic. In the noncompact case, $f$ is nonconstant by Lemma~\ref{lem:potential}. Hence $|\nabla f|^2$ is a nonzero analytic function, whose zero set has measure zero. The last assertion uses the connectedness of $M$.

We next record the coarea identities and integral scalar curvature bound used in the weighted estimate; see \cite{CZpotential,CZvolume}. Local absolute continuity allows us to integrate these identities across critical values of the potential.

\begin{lemma}\label{lem:coarea}
Let $(M^n,g,f)$ be a connected, complete, noncompact gradient shrinking Ricci soliton satisfying
\[
\Ric+\nabla^2f=\frac12g,\qquad R+|\nabla f|^2=f,
\]
where $R$ is the scalar curvature of $g$.
Set $D_t=\{x\in M:f(x)<t\}$, and let $dv_g$ and $d\sigma_g$ denote the Riemannian volume and hypersurface measures, respectively. The functions $t\mapsto\Vol(D_t)$ and $t\mapsto\int_{D_t}R\dd$ are locally absolutely continuous for $t>\min_Mf$. For every regular value $t$ in this interval,
\begin{equation}\label{eq:coareaderivatives}
\begin{aligned}
\frac d{dt}\Vol(D_t)&=\int_{\partial D_t}\frac1{|\nabla f|}\ds,\\
\frac d{dt}\int_{D_t}R\dd&=\int_{\partial D_t}\frac R{|\nabla f|}\ds.
\end{aligned}
\end{equation}
Moreover, for almost every $t>\min_Mf$,
\begin{equation}\label{eq:coareaidentity}
\frac n2\Vol(D_t)-\int_{D_t}R\dd
=t\frac d{dt}\Vol(D_t)-\frac d{dt}\int_{D_t}R\dd.
\end{equation}
For every $t>\min_Mf$ one has
\begin{equation}\label{eq:Rintegralbound}
0\le\int_{D_t}R\dd\le\frac n2\Vol(D_t).
\end{equation}
\end{lemma}
\begin{proof}
Since $f$ is proper and its critical set has measure zero by Lemma~\ref{lem:analytic}, we may apply the coarea formula on compact sublevel sets. On every compact slab, the pushforwards of $dv_g$ and $R\,dv_g$ have the integrable coarea densities in \eqref{eq:coareaderivatives}; the critical set contributes no volume measure. This proves \eqref{eq:coareaderivatives} and local absolute continuity, including across critical values.

For a regular value $t$, the divergence theorem and \eqref{eq:fidentities} yield
\[
\frac n2\Vol(D_t)-\int_{D_t}R\dd
=\int_{\partial D_t}|\nabla f|\ds.
\]
On $\partial D_t$, the normalization gives $|\nabla f|^2=t-R$. Substituting this identity and using \eqref{eq:coareaderivatives} proves \eqref{eq:coareaidentity}. The last inequality follows at regular values from the nonnegativity of the boundary integral and $R\ge0$. Continuity then extends it to every $t>\min_Mf$.
\end{proof}

\subsection{Rigidity and the number of ends}

Once we prove that $R=(n-1)/2$, the rigidity theorem of Fern\'andez-L\'opez and Garc\'ia-R\'io \cite[Theorem 10]{FLGR} gives the geometric classification: a finite quotient of an Einstein cylinder, with the following normalized potential on the product cover.

\begin{theorem}\label{thm:constant}
Let $(M^n,g,f)$ be a connected, complete, noncompact gradient shrinking Ricci soliton satisfying
\[
\Ric+\nabla^2f=\frac12g,\qquad R+|\nabla f|^2=f,
\]
where $R$ is the scalar curvature of $g$.
If
\[
R\equiv\frac{n-1}{2},
\]
then $M$ is rigid. More precisely, $n\ge3$ and $(M,g)$ is isometric to a finite quotient of $N^{n-1}\times\mathbb R$, where $N$ is compact and $\Ric_{g_N}=\tfrac12g_N$. On the product cover, the normalized potential is
\[
f(y,s)=\frac{s^2}{4}+\frac{n-1}{2}.
\]
\end{theorem}

Here rigidity is understood in the sense of Petersen--Wylie \cite{PW}; the quotient preserves the product soliton structure. In dimension two the prescribed value is $R=1/2$, which is excluded by \cite[Theorem 10]{FLGR}.

The following one-end criterion of Munteanu and Wang \cite[Theorem 4.10]{MWends}, together with Theorem~\ref{thm:main}, yields the restriction on the scalar curvature range in Corollary~\ref{cor:range}.

\begin{theorem}\label{thm:MWends}
Let $(M^n,g,f)$, $n\ge4$, be a connected, complete, noncompact gradient shrinking Ricci soliton satisfying
\[
\Ric+\nabla^2f=\frac12g,\qquad R+|\nabla f|^2=f,
\]
where $R$ is the scalar curvature of $g$.
If $R\le n/3$ on $M$, then $M$ has exactly one end.
\end{theorem}

\section{The critical scalar curvature bound}\label{sec:critical}
Throughout this section and the next, we assume the hypotheses of Theorem~\ref{thm:main} and use the normalization and notation of Section~\ref{sec:prelim}. In particular, $M$ is noncompact and $R\ge(n-1)/2$. We first determine $\delta$ and the volume growth of the sublevel sets.

Lemma~\ref{lem:volume} gives $\inf_M R\le(n-1)/2$. Together with the assumed lower bound, this yields the following proposition.

\begin{proposition}\label{p:inf}
Let $(M^n,g,f)$ be a connected, complete, noncompact gradient shrinking Ricci soliton satisfying
\[
\Ric+\nabla^2f=\frac12g,\qquad R+|\nabla f|^2=f,
\]
where $R$ is the scalar curvature of $g$.
Suppose that $R\ge(n-1)/2$ on $M$.
Writing $\delta=\inf_MR$, one has
\begin{equation}\label{eq:criticalinf}
\delta=\frac{n-1}{2}.
\end{equation}
\end{proposition}

At this endpoint, \eqref{eq:volumeupper} and \eqref{eq:volumelower} give two-sided linear volume growth. The inclusions \eqref{eq:inclusions}, together with $\sqrt t\asymp\sqrt{t-\delta}$ for large $t$, then yield the following sublevel volume estimate.

\begin{proposition}\label{p:growth}
Let $(M^n,g,f)$ be a connected, complete, noncompact gradient shrinking Ricci soliton satisfying
\[
\Ric+\nabla^2f=\frac12g,\qquad R+|\nabla f|^2=f,
\]
where $R$ is the scalar curvature of $g$.
Suppose that $R\ge(n-1)/2$ on $M$.
Set $\delta=\inf_MR$ and $D_t=\{x\in M:f(x)<t\}$. Then there exist constants $C\ge1$ and $t_*>\delta$ such that, for every $t\ge t_*$,
\begin{equation}\label{eq:sublevelgrowth}
C^{-1}\sqrt{t-\delta}\le\Vol(D_t)\le C\sqrt{t-\delta}.
\end{equation}
\end{proposition}

Set $u=R-\delta\ge0$. Then
\begin{equation}\label{eq:excess}
f-\delta=u+|\nabla f|^2\ge0,
\end{equation}
and \eqref{eq:fidentities} gives
\begin{equation}\label{eq:excesspotential}
\Delta f=\frac12-u,\qquad\Delta_f(f-\delta)=\frac12-(f-\delta).
\end{equation}

By Lemmas~\ref{lem:potential} and \ref{lem:analytic}, the potential is proper and nonconstant, and its critical set has measure zero. We may therefore apply Lemma~\ref{lem:coarea} to prove the following integrability estimate with weight $(t-\delta)^{-3/2}$. The logarithmic barrier will then yield a contradiction unless the scalar curvature excess vanishes.

\begin{proposition}\label{p:int}
Let $(M^n,g,f)$ be a connected, complete, noncompact gradient shrinking Ricci soliton satisfying
\[
\Ric+\nabla^2f=\frac12g,\qquad R+|\nabla f|^2=f,
\]
where $R$ is the scalar curvature of $g$.
Suppose that $R\ge(n-1)/2$ on $M$.
Set $\delta=(n-1)/2$, $u=R-\delta$, and $D_t=\{x\in M:f(x)<t\}$. Then, for every sufficiently large $t_0>\delta$,
\begin{equation}\label{eq:weightedint}
\int_{t_0}^\infty\left[\int_{D_t}u\dd\right](t-\delta)^{-3/2}\,dt<\infty.
\end{equation}
\end{proposition}
\begin{proof}
Applying \eqref{eq:coareaidentity} with $R=\delta+u$ and $\delta=(n-1)/2$ gives, for almost every sufficiently large $t$,
\begin{equation}\label{eq:derivativeidentity}
\frac d{dt}\left[\Vol(D_t)(t-\delta)^{-1/2}\right]
=\left[\frac d{dt}\int_{D_t}u\dd-\int_{D_t}u\dd\right](t-\delta)^{-3/2}.
\end{equation}
Integrating over $[t_0,t_1]$ and applying integration by parts to the derivative of the inner integral, we obtain
\begin{equation}\label{eq:integratedidentity}
\begin{aligned}
&\int_{t_0}^{t_1}\left[\int_{D_t}u\dd\right](t-\delta)^{-3/2}\,dt\\
&\quad=\Vol(D_{t_0})(t_0-\delta)^{-1/2}-\Vol(D_{t_1})(t_1-\delta)^{-1/2}\\
&\qquad+\left[\int_{D_{t_1}}u\dd\right](t_1-\delta)^{-3/2}
-\left[\int_{D_{t_0}}u\dd\right](t_0-\delta)^{-3/2}\\
&\qquad+\frac32\int_{t_0}^{t_1}\left[\int_{D_t}u\dd\right](t-\delta)^{-5/2}\,dt.
\end{aligned}
\end{equation}
To bound the right-hand side, we use the upper bound in \eqref{eq:sublevelgrowth}. Together with \eqref{eq:Rintegralbound} and $\delta=(n-1)/2$, it gives
\[
\begin{aligned}
0\le\int_{D_t}u\dd
&=\int_{D_t}R\dd-\delta\Vol(D_t)\\
&\le\frac12\Vol(D_t)\le C\sqrt{t-\delta}.
\end{aligned}
\]
Thus the third term on the right-hand side of \eqref{eq:integratedidentity} is bounded, while the last integral is at most
\[
C\int_{t_0}^{\infty}(t-\delta)^{-2}\,dt<\infty.
\]
Discarding the negative terms, we obtain an upper bound independent of $t_1$ for the nonnegative left-hand side. Monotone convergence then proves \eqref{eq:weightedint}.
\end{proof}

We can also prove the integrability estimate by monotonicity. Write $V(t)=\Vol(D_t)$, $U(t)=\int_{D_t}u\dd$, and $a=t-\delta$. By \eqref{eq:coareaidentity}, $aV'=V/2-U+U'$ almost everywhere. Hence the locally absolutely continuous function
\[
A(t)=V(t)a^{-1/2}-U(t)a^{-3/2}
\]
satisfies
\[
A'(t)=-U(t)a^{-3/2}\left(1-\frac{3}{2a}\right).
\]
Since $0\le U\le V/2$, for $a\ge3$ we have $A\ge0$ and $A'\le-\tfrac12Ua^{-3/2}$. Integration gives
\[
\int_{t_0}^T U(t)(t-\delta)^{-3/2}\,dt\le2A(t_0)
\]
whenever $t_0-\delta\ge3$. This proves the same integrability estimate without using the upper volume bound.

We next derive a sum-of-squares identity for the scalar curvature excess on the regular set of the potential. Its differential inequality will be used in the maximum principle and barrier arguments.

\begin{lemma}\label{lem:differential}
Let $(M^n,g,f)$ be a connected, complete, noncompact gradient shrinking Ricci soliton satisfying
\[
\Ric+\nabla^2f=\frac12g,\qquad R+|\nabla f|^2=f,
\]
where $R$ is the scalar curvature of $g$.
Suppose that $R\ge(n-1)/2$ on $M$.
Set $\delta=(n-1)/2$, $u=R-\delta$, and $\Delta_f=\Delta-\langle\nabla f,\nabla(\cdot)\rangle$. At any point of $\{|\nabla f|>0\}$, choose an orthonormal frame with $e_1=\nabla f/|\nabla f|$. Then
\begin{equation}\label{eq:squareidentity}
\begin{aligned}
&\left(\Delta_f-\frac{\langle\nabla f,\nabla(\cdot)\rangle}{|\nabla f|^2}\right)u\\
&\quad=-u-\frac2{n-1}\left[u-\Ric(e_1,e_1)\right]^2-2\Ric(e_1,e_1)^2\\
&\qquad-4\sum_{i=2}^n\Ric(e_1,e_i)^2\\
&\qquad-2\sum_{i,j=2}^n\left[\Ric(e_i,e_j)-\left(\frac12+\frac{u-\Ric(e_1,e_1)}{n-1}\right)g(e_i,e_j)\right]^2.
\end{aligned}
\end{equation}
In particular,
\begin{equation}\label{eq:differentialineq}
\left(\Delta_f-\frac{\langle\nabla f,\nabla(\cdot)\rangle}{|\nabla f|^2}\right)u\le-u.
\end{equation}
Moreover,
\begin{equation}\label{eq:modifiedpotential}
\left(\Delta_f-\frac{\langle\nabla f,\nabla(\cdot)\rangle}{|\nabla f|^2}\right)(f-\delta)=-(f-\delta)-\frac12.
\end{equation}
\end{lemma}
\begin{proof}
From \eqref{eq:Ridentities},
\[
\Delta_fu=\delta+u-2|\Ric|^2,\qquad
\frac{\langle\nabla f,\nabla u\rangle}{|\nabla f|^2}=2\Ric(e_1,e_1).
\]
Decomposing the Ricci tensor in the chosen frame, we have
\[
\begin{aligned}
|\Ric|^2&=\Ric(e_1,e_1)^2+2\sum_{i=2}^n\Ric(e_1,e_i)^2\\
&\quad+\sum_{i,j=2}^n\left[\Ric(e_i,e_j)-\left(\frac12+\frac{u-\Ric(e_1,e_1)}{n-1}\right)g(e_i,e_j)\right]^2\\
&\quad+\frac{[\delta+u-\Ric(e_1,e_1)]^2}{n-1}.
\end{aligned}
\]
Since $\delta=(n-1)/2$, the tangential average is
\[
\frac{\delta+u-\Ric(e_1,e_1)}{n-1}
=\frac12+\frac{u-\Ric(e_1,e_1)}{n-1}.
\]
Substitution proves \eqref{eq:squareidentity}. Finally, \eqref{eq:excesspotential} and $\langle\nabla f,\nabla(f-\delta)\rangle/|\nabla f|^2=1$ give \eqref{eq:modifiedpotential}.
\end{proof}
\begin{remark}\label{rem:singularoperator}
The coefficients of this operator are smooth on $\{|\nabla f|>0\}$. At a critical point, $u=f-\delta$ by \eqref{eq:excess}; this identity will exclude critical points from the minimum argument in Lemma~\ref{lem:barrier}.
\end{remark}

\section{The logarithmic barrier and the Liouville theorem}\label{sec:barrier}

We first use the strong minimum principle and analyticity to show that the scalar curvature excess either vanishes identically or is strictly positive on the exterior region.

\begin{lemma}\label{lem:positive}
Let $(M^n,g,f)$ be a connected, complete, noncompact gradient shrinking Ricci soliton satisfying
\[
\Ric+\nabla^2f=\frac12g,\qquad R+|\nabla f|^2=f,
\]
where $R$ is the scalar curvature of $g$.
Suppose that $R\ge(n-1)/2$ on $M$.
Set $\delta=(n-1)/2$ and $u=R-\delta\ge0$. Then either $u\equiv0$ on $M$, or
\begin{equation}\label{eq:positivity}
u>0\qquad\text{on }\{f>\delta\}.
\end{equation}
\end{lemma}
\begin{proof}
Suppose that $u$ vanishes at a point where $f>\delta$. By \eqref{eq:excess}, $|\nabla f|^2=f-\delta>0$ at that point. The differential expression in \eqref{eq:differentialineq} is therefore smooth and uniformly elliptic in a neighborhood. The strong minimum principle implies that $u$ vanishes on an open set. Since $u=R-\delta$ is analytic, Lemma~\ref{lem:analytic} yields $u\equiv0$ on $M$.
\end{proof}

In the positive case, we use the differential inequality to obtain a logarithmic lower bound at infinity. Together with sublevel volume growth, this will contradict the weighted integrability estimate.

\begin{lemma}\label{lem:barrier}
Let $(M^n,g,f)$ be a connected, complete, noncompact gradient shrinking Ricci soliton satisfying
\[
\Ric+\nabla^2f=\frac12g,\qquad R+|\nabla f|^2=f,
\]
where $R$ is the scalar curvature of $g$.
Suppose that $R\ge(n-1)/2$ on $M$.
Set $\delta=(n-1)/2$ and $u=R-\delta\ge0$, and suppose that $u\not\equiv0$ on $M$. Then there exist a sufficiently large regular value $t_0>\delta+e$ of $f$ and a constant $\varepsilon>0$ such that
\begin{equation}\label{eq:logbarrier}
u(x)\ge\frac{\varepsilon}{\log(f(x)-\delta)}\qquad\text{whenever }f(x)\ge t_0.
\end{equation}
\end{lemma}
\begin{proof}
The proper nonconstant function $f$ is unbounded, and $f(M)$ is an interval. By Sard's theorem, we may choose a sufficiently large regular value $t_0>\delta+e$ in the interior of $f(M)$ so that Proposition~\ref{p:int} applies with this lower limit. The corresponding level set is nonempty and compact. By Lemma~\ref{lem:positive},
\[
\min_{\{f=t_0\}}u>0.
\]
Choose
\begin{equation}\label{eq:epsilonchoice}
0<\varepsilon<\min\left\{\log(t_0-\delta)\min_{\{f=t_0\}}u,\ (t_0-\delta)\log(t_0-\delta)\right\}.
\end{equation}
Then
\[
u-\frac{\varepsilon}{\log(f-\delta)}
\]
is nonnegative on $\{f=t_0\}$ and has a nonnegative lower limit at infinity. If it were negative somewhere on $\{f\ge t_0\}$, properness of $f$ and this lower limit would imply that it attains a negative global minimum in a compact sublevel set, away from the boundary level. This minimum cannot lie in $\{|\nabla f|=0\}$: at a critical point, \eqref{eq:excess} gives $u=f-\delta$, while \eqref{eq:epsilonchoice} gives
\[
(f-\delta)-\frac{\varepsilon}{\log(f-\delta)}>0.
\]
Here we use that $s\mapsto s\log s$ is increasing for $s>e$ and that $f-\delta\ge t_0-\delta>e$. Thus a negative minimum lies in $\{|\nabla f|>0\}$. At this point, all first-order terms vanish and the Hessian is nonnegative, so
\begin{equation}\label{eq:minimumineq}
\left(\Delta_f-\frac{\langle\nabla f,\nabla(\cdot)\rangle}{|\nabla f|^2}\right)
\left[u-\frac{\varepsilon}{\log(f-\delta)}\right]\ge0.
\end{equation}
On $\{f-\delta>1\}$, direct differentiation and \eqref{eq:modifiedpotential} give
\begin{equation}\label{eq:logderivative}
\begin{aligned}
&\left(\Delta_f-\frac{\langle\nabla f,\nabla(\cdot)\rangle}{|\nabla f|^2}\right)
\left[\frac1{\log(f-\delta)}\right]\\
&\quad=\frac{f-\delta+\tfrac12}{(f-\delta)[\log(f-\delta)]^2}\\
&\qquad+\left[\frac1{(f-\delta)^2[\log(f-\delta)]^2}
+\frac2{(f-\delta)^2[\log(f-\delta)]^3}\right]|\nabla f|^2>0.
\end{aligned}
\end{equation}
By \eqref{eq:differentialineq} and \eqref{eq:logderivative}, the left-hand side of \eqref{eq:minimumineq} is strictly negative, a contradiction.
\end{proof}

\begin{proof}[Proof of Theorem~\ref{thm:main}]
Suppose, for a contradiction, that $u\not\equiv0$. For all sufficiently large $t>t_0$, the logarithmic barrier gives
\begin{equation}\label{eq:shelllower}
\begin{aligned}
\int_{D_t}u\dd
&\ge\int_{\{t_0\le f<t\}}\frac{\varepsilon}{\log(f-\delta)}\dd\\
&\ge\frac{\varepsilon}{\log(t-\delta)}\left[\Vol(D_t)-\Vol(D_{t_0})\right].
\end{aligned}
\end{equation}
Using the lower bound in \eqref{eq:sublevelgrowth}, we obtain
\begin{equation}\label{eq:finallower}
\int_{D_t}u\dd\ge C\frac{\sqrt{t-\delta}}{\log(t-\delta)}
\end{equation}
for all sufficiently large $t$. Therefore
\[
\begin{aligned}
\int_{t_0}^\infty\left[\int_{D_t}u\dd\right](t-\delta)^{-3/2}\,dt
&\ge C\int^{\infty}\frac{dt}{(t-\delta)\log(t-\delta)}=\infty.
\end{aligned}
\]
This contradicts Proposition~\ref{p:int}. Therefore $u\equiv0$, or equivalently $R\equiv(n-1)/2$. Theorem~\ref{thm:constant} now gives the geometric conclusion.
\end{proof}


\section{Critical-point criteria for one-endedness}\label{sec:criticalends}

In this section and the next, we drop the scalar lower bound used in Sections~3 and~4. Unless otherwise stated, $(M^n,g,f)$ is a connected, complete, noncompact gradient shrinking Ricci soliton satisfying \eqref{eq:normalization}. We first establish a topological criterion involving only the critical points of $f$, including degenerate ones.

\subsection{Ends and a compact core}
For a connected noncompact manifold, the ends are defined by the inverse system of unbounded connected components of complements of compact sets. In particular, one-endedness means that sufficiently large compact sets cannot separate two distinct directions to infinity. In the topological arguments, homology and cohomology are singular homology and cohomology with coefficients in $\mathbb Z_2$, unless otherwise indicated. With these coefficients, no orientability assumption is needed.

The next lemma shows that two strictly positive Hessian directions at each critical point suffice for one-endedness. The proof uses a compact core and applies also to degenerate critical points.

\begin{lemma}\label{lem:compactcore}
Let $X$ be a connected, noncompact, smooth $n$-manifold without boundary, equipped with an auxiliary smooth Riemannian metric. Let $h:X\to\mathbb R$ be a smooth proper function that is bounded below. Suppose that its critical set $\operatorname{Crit}(h)=\{x\in X:dh_x=0\}$ is compact and that, at every point of this set, its Hessian has at least two strictly positive eigenvalues. Then $X$ has exactly one end. Moreover, $H_{n-1}(X;\mathbb Z_2)=0$.
\end{lemma}
\begin{proof}
With respect to the given auxiliary metric, let $C=\operatorname{Crit}(h)$. Choose a regular value $a$ larger than all critical values, and put $K=\{h\le a\}$. Then $K$ is a compact smooth manifold with boundary. On $\{h\ge a\}$ the vector field
\[
Y=\frac{\nabla h}{|\nabla h|^2}
\]
satisfies $Y(h)=1$. Its flow exists across every finite interval of $h$-values because such a slab is compact and has no critical points. Consequently,
\begin{equation}\label{eq:collar}
X\setminus\operatorname{int}K\cong\partial K\times[0,\infty).
\end{equation}
Collapsing the second factor gives a deformation retraction of $X$ onto $K$. Thus $K$ is connected and the number of ends of $X$ equals the number of connected components of $\partial K$.

To apply Morse theory, we perturb $h$ while preserving two positive Hessian directions near its critical set. At each point of the compact set $C$, the second largest Hessian eigenvalue is positive. By continuity and compactness, there is a relatively compact neighborhood $U$ of $C$, with closure in $\operatorname{int}K$, on which the Hessian has two positive eigenvalues bounded away from zero. After shrinking $U$ if necessary, a sufficiently small uniform $C^2$ perturbation preserves these two positive eigenvalues on $U$. On the compact set $K\setminus U$ the gradient of $h$ has a positive lower bound. A sufficiently small $C^1$ perturbation therefore creates no critical points there.

By relative Morse approximation, obtained from parametric transversality applied to the differential section, we may choose a Morse function $H$ arbitrarily $C^2$ close to $h$, with $H-h$ supported in a fixed compact set $S\Subset\operatorname{int}K$; see \cite[Chapter~3, Theorem~2.7, and Chapter~6, Theorem~1.2]{Hirsch}. The preceding gradient and Hessian bounds ensure that all critical points of $H$ lie in $U$ and retain two positive Hessian directions. Since $a-\max_S h>0$, taking the perturbation sufficiently small also gives $H<a$ on $S$. Together with $H=h$ off $S$, this yields $\{H\le a\}=K$.

At every critical point of $H$ there are at least two positive Hessian eigenvalues. Thus
\begin{equation}\label{eq:indexbound}
\operatorname{ind}_H(p)\le n-2.
\end{equation}
Starting below the minimum of $H$, the cell-attachment theorem of Morse theory gives $K$ the homotopy type of a finite CW complex with one cell for each critical point, of dimension equal to its Morse index; see \cite[Section~3]{Milnor}. By \eqref{eq:indexbound}, this complex has dimension at most $n-2$, and hence
\begin{equation}\label{eq:homologyvanish}
H_{n-1}(K;\mathbb Z_2)=0.
\end{equation}

Write $b=\#\pi_0(\partial K)$. Since $K$ is compact and connected with nonempty boundary, mod-two Poincar\'e--Lefschetz duality and the long exact sequence of the pair give
\[
\underbrace{H_n(K,\partial K;\mathbb Z_2)}_{\cong\mathbb Z_2}
\longrightarrow
\underbrace{H_{n-1}(\partial K;\mathbb Z_2)}_{\cong\mathbb Z_2^{\,b}}
\longrightarrow H_{n-1}(K;\mathbb Z_2)=0.
\]
The first map must be surjective, so $b\le1$. Noncompactness and \eqref{eq:collar} imply $b\ge1$. Thus $b=1$, and $X$ has exactly one end. The deformation retraction and \eqref{eq:homologyvanish} also give the claimed homology vanishing.
\end{proof}

\begin{remark}\label{rem:twopositive}
The strict positivity of the two Hessian eigenvalues is essential. In the cylinder example, the Hessian has only one positive eigenvalue and $n-1$ zero eigenvalues along the critical hypersurface, so two nonnegative eigenvalues do not suffice.
\end{remark}

\subsection{Ricci nonnegativity restricted to the critical set}
We now prove the critical-set theorem announced in the introduction.

\begin{proof}[Proof of Theorem~\ref{thm:criticalintro}]
Put $c=(n-1)/2$. On the critical set the normalization gives
\begin{equation}\label{eq:criticalnormalization}
f(p)=R(p),\qquad \nabla^2f(p)=\tfrac12g_p-\operatorname{Ric}_p.
\end{equation}
Consequently $\operatorname{Crit}(f)\subset\{f<c\}$. Since $f$ is proper, the closed set $\operatorname{Crit}(f)$ is contained in the compact set $\{f\le c\}$ and is therefore compact. Thus a uniform gap on the critical set follows from the pointwise strict inequality and compactness.

At $p\in\operatorname{Crit}(f)$ write the Ricci eigenvalues in increasing order,
\[
0\le\lambda_1(p)\le\lambda_2(p)\le\cdots\le\lambda_n(p).
\]
If $\lambda_2(p)\ge1/2$, then
\[
R(p)=\lambda_1(p)+\sum_{i=2}^n\lambda_i(p)\ge\frac{n-1}{2},
\]
contrary to the hypothesis. Thus $\lambda_1(p),\lambda_2(p)<1/2$, and \eqref{eq:criticalnormalization} shows that $\nabla^2f(p)$ has at least two strictly positive eigenvalues. Applying Lemma~\ref{lem:compactcore} to $h=f$ gives one-endedness and $H_{n-1}(M;\mathbb Z_2)=0$.
\end{proof}

Applying the critical-set theorem under a global scalar upper bound gives the following corollary, with Ricci nonnegativity required only at critical points.

\begin{corollary}\label{cor:globalcritical}
Let $(M^n,g,f)$ be a connected, complete, noncompact gradient shrinking Ricci soliton satisfying
\[
\Ric+\nabla^2f=\frac12g,\qquad R+|\nabla f|^2=f,
\]
where $R$ is the scalar curvature of $g$.
Suppose that $R<(n-1)/2$ on $M$ and that $\operatorname{Ric}_p\ge0$ for every $p$ in the critical set $\operatorname{Crit}(f)=\{x\in M:df_x=0\}$. Then $M$ has exactly one end.
\end{corollary}

\begin{remark}\label{rem:localRic}
In Theorem~\ref{thm:criticalintro}, Ricci nonnegativity is required only on the critical set. The result applies without bounded Riemann curvature, nondegeneracy of the critical set, a K\"ahler structure, or the scalar lower bound of Theorem~\ref{thm:main}.
\end{remark}

\section{Ricci norm pinching and the cylindrical alternative}\label{sec:norm}

We now consider the non-strict global Ricci norm bound. The global hypothesis allows us to use the scalar curvature equation to prove rigidity in the equality case.

\subsection{A scalar upper bound from the Ricci norm}

We first derive the scalar upper bound needed for the critical-point argument and show that equality at one point forces constant scalar curvature.

\begin{lemma}\label{lem:scalarupper}
Let $(M^n,g,f)$ be a connected, complete, noncompact gradient shrinking Ricci soliton satisfying
\[
\Ric+\nabla^2f=\frac12g,\qquad R+|\nabla f|^2=f,
\]
where $R$ is the scalar curvature of $g$.
Suppose that
\[
|\Ric|^2\le\frac{n-1}{4}\qquad\text{on }M,
\]
and put $c=(n-1)/2$. Then $R\le c$ on $M$. Either $R<c$ everywhere, or $R\equiv c$ everywhere.
\end{lemma}
\begin{proof}
By Cauchy--Schwarz and scalar nonnegativity,
\begin{equation}\label{eq:roughscalar}
0\le R\le\sqrt n\,|\operatorname{Ric}|\le\frac{\sqrt{n(n-1)}}2.
\end{equation}
In particular $R$ is bounded, whereas $f$ is proper and tends to infinity. The scalar equation gives
\begin{equation}\label{eq:scalardifferential}
\Delta_fR=R-2|\operatorname{Ric}|^2\ge R-c.
\end{equation}
For $\varepsilon>0$, the smooth function $G_\varepsilon=R-\varepsilon f$ tends to $-\infty$ at infinity, and therefore attains its maximum at a point $p_\varepsilon$. At that point its gradient vanishes and its Laplacian is nonpositive. Using $\Delta_ff=n/2-f$, we obtain
\[
0\ge\Delta_fG_\varepsilon(p_\varepsilon)
\ge R(p_\varepsilon)-c-\varepsilon\bigl(n/2-f(p_\varepsilon)\bigr),
\]
or
\begin{equation}\label{eq:maxpoint}
R(p_\varepsilon)+\varepsilon f(p_\varepsilon)
\le c+\varepsilon n/2.
\end{equation}
Since $f=R+|\nabla f|^2\ge0$, for every fixed $x\in M$ it follows that
\[
R(x)-\varepsilon f(x)\le G_\varepsilon(p_\varepsilon)
\le c+\varepsilon n/2.
\]
Letting $\varepsilon\downarrow0$ gives $R(x)\le c$, without assuming attainment of the scalar supremum or bounds on curvature derivatives.

If equality holds at some point, let $v=c-R\ge0$. Equation~\eqref{eq:scalardifferential} says
\[
(\Delta_f-1)v\le0.
\]
Since the zeroth-order coefficient is $-1$ and $M$ is connected, the strong minimum principle gives $v\equiv0$. Otherwise $v>0$ everywhere.
\end{proof}

\subsection{The strict scalar case}

When the scalar curvature is strictly below the cylindrical value, the critical set is compact and the norm bound forces two strictly positive Hessian directions at every critical point. Lemma~\ref{lem:compactcore} then gives one-endedness and the homology vanishing.

\begin{lemma}\label{lem:normtwopositive}
Let $(M^n,g,f)$ be a connected, complete, noncompact gradient shrinking Ricci soliton satisfying
\[
\Ric+\nabla^2f=\frac12g,\qquad R+|\nabla f|^2=f,
\]
where $R$ is the scalar curvature of $g$.
Suppose that
\[
|\Ric|^2\le\frac{n-1}{4},\qquad R<\frac{n-1}{2}
\qquad\text{everywhere on }M.
\]
Then $\nabla^2f$ has at least two strictly positive eigenvalues at every critical point of $f$. Consequently, $M$ has exactly one end and $H_{n-1}(M;\mathbb Z_2)=0$.
\end{lemma}
\begin{proof}
Suppose that $\lambda_2\ge1/2$ at a critical point, with the Ricci eigenvalues ordered increasingly. Then
\[
\frac{n-1}{4}\ge\sum_{i=1}^n\lambda_i^2
\ge\lambda_1^2+\frac{n-1}{4}.
\]
Equality holds throughout, so $\lambda_1=0$ and $\lambda_2=\cdots=\lambda_n=1/2$. Thus $R=(n-1)/2$ at this point, contrary to the strict scalar inequality. Hence $\lambda_2<1/2$ at every critical point, and $\nabla^2f=\tfrac12g-\Ric$ has at least two strictly positive eigenvalues there. Put $c=(n-1)/2$. Since $f=R<c$ on $\Crit(f)$ and $f$ is proper, the closed set $\Crit(f)$ is contained in the compact set $\{f\le c\}$ and is therefore compact. Lemma~\ref{lem:compactcore} now gives the result.
\end{proof}

\subsection{The equality case and the ends of a cylinder quotient}

In the constant-scalar case, Theorem~\ref{thm:constant} identifies the soliton as a finite quotient of an Einstein cylinder. We determine the number of ends of the quotient and the condition under which it splits as a global product.

\begin{lemma}\label{lem:cylinderquotient}
Let $(M^n,g,f)$ be a connected, complete, noncompact gradient shrinking Ricci soliton satisfying
\[
\Ric+\nabla^2f=\frac12g,\qquad R+|\nabla f|^2=f,
\]
where $R$ is the scalar curvature of $g$.
Suppose that
\[
R\equiv\frac{n-1}{2}\qquad\text{on }M.
\]
Then $n\ge3$, and $(M,g)$ is isometric to a finite quotient of $(N^{n-1}\times\mathbb R,g_N+ds^2)$, where $(N,g_N)$ is a connected compact Einstein manifold satisfying $\Ric_{g_N}=\tfrac12g_N$. In particular, $M$ has at most two ends. If $M$ has more than one end, then it is itself a global product $(\widehat N^{n-1}\times\mathbb R,g_{\widehat N}+ds^2)$, with $\widehat N$ compact and $\Ric_{g_{\widehat N}}=\tfrac12g_{\widehat N}$, and has exactly two ends.
\end{lemma}
\begin{proof}
Theorem~\ref{thm:constant} implies $n\ge3$ and rigidity. On the universal cover we may write
\[
(\widetilde M,\widetilde g,\widetilde f)
=\left(N\times\mathbb R,g_N+ds^2,\frac{s^2}{4}+\frac{n-1}{2}\right),
\]
after translating the line coordinate. The Einstein factor $N$ is complete and satisfies $\operatorname{Ric}_N=\tfrac12g_N$, so it is compact by Myers' theorem. For any deck transformation $\gamma$, the lifted potential is invariant. The product soliton structure, equivalently the kernel and image of the Ricci tensor on this cover, is also invariant. Hence
\begin{equation}\label{eq:deckaction}
\gamma(y,s)=(\phi_\gamma(y),\epsilon_\gamma s),\qquad
\epsilon_\gamma\in\{1,-1\},
\end{equation}
where $\phi_\gamma$ is an isometry of $N$. Indeed, a line isometry has the form $s\mapsto\epsilon s+b$, and invariance of $s^2/4+(n-1)/2$ for all $s$ forces $b=0$.

Since the deck group acts properly discontinuously and preserves the compact set $N\times\{0\}$, it is finite. Let $\Gamma_+$ be the subgroup with $\epsilon_\gamma=1$. It acts freely on $N$, and
\[
(N\times\mathbb R)/\Gamma_+=(N/\Gamma_+)\times\mathbb R.
\]
If the deck group equals $\Gamma_+$, the quotient is the asserted global product and has exactly two ends. Otherwise the residual group has order two and exchanges the two ends of this product. For every sufficiently large $A$, the regions $s>A$ and $s<-A$ are identified in the full quotient, so its exterior is connected. Thus the full quotient has exactly one end. This proves the lemma.
\end{proof}

\begin{proof}[Proof of Theorem~\ref{thm:normintro}]
Lemma~\ref{lem:scalarupper} gives either $R<c$ everywhere or $R\equiv c$, where $c=(n-1)/2$. The first case is one-ended by Lemma~\ref{lem:normtwopositive}. In the second case Lemma~\ref{lem:cylinderquotient} gives either one end or a global Einstein cylinder with exactly two ends. If
\[
|\operatorname{Ric}|^2(x_0)<\frac{n-1}{4}
\]
at some point $x_0$, the scalar equation excludes the constant-scalar case, in which $|\operatorname{Ric}|^2\equiv R/2=(n-1)/4$. This proves the final assertion.
\end{proof}

Since the cylindrical alternative in Theorem~\ref{thm:normintro} attains the norm threshold everywhere, a pointwise strict bound gives the following one-end criterion.

\begin{corollary}\label{cor:strictglobalnorm}
Let $(M^n,g,f)$ be a connected, complete, noncompact gradient shrinking Ricci soliton satisfying
\[
\Ric+\nabla^2f=\frac12g,\qquad R+|\nabla f|^2=f,
\]
where $R$ is the scalar curvature of $g$.
If
\[
|\Ric|(x)<\frac{\sqrt{n-1}}{2}\qquad\text{for every }x\in M,
\]
then $M$ has exactly one end. A uniform strict gap is not required.
\end{corollary}

\subsection{Sharpness and comparison with the literature}
For $n\ge3$, let $(N^{n-1},g_N)$ be any connected compact Einstein manifold with $\operatorname{Ric}_{g_N}=\tfrac12g_N$. On its product with $\mathbb R$,
\begin{equation}\label{eq:cylinderexample}
R\equiv\frac{n-1}{2},\quad
|\operatorname{Ric}|^2\equiv\frac{n-1}{4},\quad
f(y,s)=\frac{s^2}{4}+\frac{n-1}{2}.
\end{equation}
This complete smooth example has two ends and attains the norm threshold, showing that the cylindrical alternative is necessary under the non-strict global bound. Its critical set is $N\times\{0\}$, where the Hessian has one eigenvalue $1/2$ and $n-1$ zero eigenvalues. The example also shows why the scalar bound in Theorem~\ref{thm:criticalintro} must be strict.

The strict global bound in Corollary~\ref{cor:strictglobalnorm} improves the sufficient numerical bound $|\operatorname{Ric}|<(n-2)/(2\sqrt n)$ in Qu--Wu \cite{QW}. The proof here is based on the critical-point Hessian, rather than on a global weighted harmonic-function estimate. The K\"ahler one-end theorem of Munteanu--Wang \cite{MWkahler} has no such Ricci norm restriction; its complex geometric setting provides a different mechanism for one-endedness.

Bertellotti--Buzano \cite{BBends} study quantitative bounds and stability of the number of ends, including singular limits and asymptotically conical ends. Their hypotheses excluding critical points sufficiently far out are automatically satisfied under the critical-set scalar hypothesis of Theorem~\ref{thm:criticalintro}. The index restriction in our proof further shows that the exterior level set is connected. Their line-splitting results for limits along suitable trajectories in \cite{BBstructure} concern the geometry at infinity; the global splitting in Theorem~\ref{thm:normintro} instead comes from the scalar equality case and the deck-group analysis above.

\begin{remark}\label{rem:noncompactautomatic}
In Theorems~\ref{thm:criticalintro} and~\ref{thm:normintro}, noncompactness follows from the curvature hypotheses. If $M$ were compact and $|\operatorname{Ric}|^2\le(n-1)/4$, then \eqref{eq:roughscalar} would imply $R<n/2$ everywhere, contrary to the integral of $\Delta f=n/2-R$. Under the critical-set hypotheses of Theorem~\ref{thm:criticalintro}, a maximum point of $f$ on a compact $M$ would satisfy $\nabla^2f\le0$, hence $\operatorname{Ric}\ge\tfrac12g$ and $R\ge n/2$, again a contradiction.
\end{remark}

\subsection{Lower Ricci norm bounds and scalar comparison}\label{subsec:lowernorm}
In this subsection and the next, $M$ is connected and complete and may be compact unless otherwise stated. We retain the normalization \eqref{eq:normalization}, but drop the standing noncompactness assumption from Sections~\ref{sec:criticalends} and~\ref{sec:norm}.

We use the scalar equation to derive a scalar lower bound from a global lower bound for $|\Ric|^2$. This comparison is the main estimate for the rigidity theorem below and requires neither a curvature upper bound nor attainment of the scalar infimum.

\begin{lemma}\label{lem:normlowercomparison}
Let $(M^n,g,f)$ be a connected, complete gradient shrinking Ricci soliton satisfying
\[
\Ric+\nabla^2f=\frac12g,\qquad R+|\nabla f|^2=f,
\]
where $R$ is the scalar curvature of $g$.
Let $a\ge0$ be a constant, and suppose that $|\Ric|^2\ge a$ everywhere on $M$. Then
\begin{equation}\label{eq:normlowerscalar}
R\ge2a\qquad\text{on }M.
\end{equation}
Either $R>2a$ everywhere, or $R\equiv2a$ and $|\Ric|^2\equiv a$ everywhere. In particular, if $|\Ric|^2>a$ at even one point, then $R>2a$ at every point.
\end{lemma}
\begin{proof}
The scalar equation and the normalization give
\begin{equation}\label{eq:lowerPDE}
\Delta_fR=R-2|\Ric|^2\le R-2a,
\qquad\Delta_ff=\frac n2-f.
\end{equation}
Scalar nonnegativity gives $R\ge0$ and hence $f=R+|\nabla f|^2\ge0$. Fix $\varepsilon>0$ and consider
\[
H_\varepsilon=R+\varepsilon f.
\]
The function $H_\varepsilon$ attains its minimum when $M$ is compact. When $M$ is noncompact, Lemma~\ref{lem:potential} gives a proper exhaustion $f$; since $R\ge0$, we have $H_\varepsilon\to+\infty$ at infinity. Thus $H_\varepsilon$ attains its minimum in either case. Let $p_\varepsilon$ be a global minimum point.

At $p_\varepsilon$ the gradient of $H_\varepsilon$ vanishes and its Laplacian is nonnegative. Equation~\eqref{eq:lowerPDE} therefore yields
\[
0\le\Delta_fH_\varepsilon(p_\varepsilon)
\le R(p_\varepsilon)-2a
+\varepsilon\left(\frac n2-f(p_\varepsilon)\right).
\]
Rearranging gives
\begin{equation}\label{eq:loweratmin}
R(p_\varepsilon)\ge2a-\frac{\varepsilon n}{2}
+\varepsilon f(p_\varepsilon).
\end{equation}
For every fixed $x\in M$, the minimizing property and $f\ge0$ imply
\begin{align*}
R(x)+\varepsilon f(x)
&\ge R(p_\varepsilon)+\varepsilon f(p_\varepsilon)\\
&\ge2a-\frac{\varepsilon n}{2}+2\varepsilon f(p_\varepsilon)\\
&\ge2a-\frac{\varepsilon n}{2}.
\end{align*}
Letting $\varepsilon\downarrow0$ with $x$ fixed proves \eqref{eq:normlowerscalar}, without any upper bound for $R$, $|\Ric|^2$ or the Riemann curvature tensor.

To prove the alternative, put $u=R-2a\ge0$. Then
\begin{equation}\label{eq:lowerequalityoperator}
(\Delta_f-1)u=-2(|\Ric|^2-a)\le0.
\end{equation}
If $u$ vanishes at some point, the strong minimum principle gives $u\equiv0$, since the zeroth-order coefficient is $-1$ and $M$ is connected. The scalar equation then yields $|\Ric|^2\equiv a$. Otherwise $u>0$ everywhere, which also proves the final assertion.
\end{proof}

\begin{proof}[Proof of Theorem~\ref{thm:normlowerintro}]
Take $a=(n-1)/4$ in Lemma~\ref{lem:normlowercomparison}. This gives $R\ge(n-1)/2$ on $M$. If $M$ is noncompact, Theorem~\ref{thm:main} gives $R\equiv(n-1)/2$, and \eqref{eq:Ridentities} then implies
\begin{equation}\label{eq:normlowerrigidequalities}
|\Ric|^2\equiv\frac{n-1}{4}.
\end{equation}
The geometric conclusion follows from Theorem~\ref{thm:constant}. In particular, this noncompact case requires $n\ge3$.

If the assumed global norm lower bound is strict at even one point, the equality alternative in Lemma~\ref{lem:normlowercomparison} is impossible. Thus $R>(n-1)/2$ everywhere, and Corollary~\ref{cor:compact} implies compactness. Equivalently, strictness at one point contradicts \eqref{eq:normlowerrigidequalities} in the noncompact case.
\end{proof}

Under the global lower bound in Theorem~\ref{thm:normlowerintro}, strictness at one point therefore gives the following compactness criterion.

\begin{corollary}\label{cor:lowerstrictcompact}
Let $(M^n,g,f)$ be a connected, complete gradient shrinking Ricci soliton satisfying
\[
\Ric+\nabla^2f=\frac12g,\qquad R+|\nabla f|^2=f,
\]
where $R$ is the scalar curvature of $g$.
Suppose that
\[
|\Ric|^2\ge\frac{n-1}{4}\qquad\text{on }M
\]
and $|\Ric|^2(x_0)>(n-1)/4$ for at least one point $x_0\in M$. Then $M$ is compact. In particular, the pointwise strict inequality $|\Ric|^2>(n-1)/4$ on $M$ implies compactness, without a uniform positive gap.
\end{corollary}

\begin{remark}\label{rem:lowerends}
In the noncompact lower-bound case, Lemma~\ref{lem:cylinderquotient} shows that $M$ has one or two ends, with two ends precisely when the cylinder quotient is a global product with a line. In particular, a simply connected $M$ is a global Einstein cylinder and has exactly two ends. The rigidity conclusion also gives $\nabla\Ric=0$ and nonnegative Ricci curvature, although neither is assumed.
\end{remark}

Applying Lemma~\ref{lem:normlowercomparison} with $a=\inf_M|\Ric|^2$ gives $2\inf_M|\Ric|^2\le\inf_MR$. Combining this with Lemmas~\ref{lem:scalar} and~\ref{lem:volume}, and using Theorem~\ref{thm:normlowerintro} for the equality case, we obtain the following corollary.

\begin{corollary}\label{cor:norminfimum}
Let $(M^n,g,f)$ be a connected, complete gradient shrinking Ricci soliton satisfying
\[
\Ric+\nabla^2f=\frac12g,\qquad R+|\nabla f|^2=f,
\]
where $R$ is the scalar curvature of $g$.
Then
\begin{equation}\label{eq:generalnorminfimum}
2\inf_M|\Ric|^2\le\inf_MR\le\frac n2,
\qquad\inf_M|\Ric|^2\le\frac n4.
\end{equation}
If $M$ is noncompact, then
\begin{equation}\label{eq:noncompactnorminfimum}
\inf_M|\Ric|^2\le\frac{n-1}{4}.
\end{equation}
In the noncompact case, equality in \eqref{eq:noncompactnorminfimum} holds if and only if $n\ge3$ and $(M,g)$ is isometric to a finite quotient of $(N^{n-1}\times\mathbb R,g_N+ds^2)$, where $(N,g_N)$ is a connected compact Einstein manifold satisfying $\Ric_{g_N}=\tfrac12g_N$.
\end{corollary}

At the larger threshold $n/4$, the same comparison forces the scalar infimum to reach its maximal value. The rigidity result of Pigola--Rimoldi--Setti and Mai--Ou \cite{PRS,MO}, recorded in Lemma~\ref{lem:scalar}, then yields a compact Einstein shrinker.

\begin{corollary}\label{cor:Einsteinnormendpoint}
Let $(M^n,g,f)$ be a connected, complete gradient shrinking Ricci soliton satisfying
\[
\Ric+\nabla^2f=\frac12g,\qquad R+|\nabla f|^2=f,
\]
where $R$ is the scalar curvature of $g$.
If
\begin{equation}\label{eq:Einsteinnormpinch}
|\Ric|^2\ge\frac n4\qquad\text{on }M,
\end{equation}
then $M$ is compact, $f$ is constant and $\Ric=\tfrac12g$. Consequently, $|\Ric|^2\equiv n/4$; strictness anywhere in \eqref{eq:Einsteinnormpinch} is impossible. Equivalently, equality $\inf_M|\Ric|^2=n/4$ characterizes compact Einstein shrinkers among the connected, complete gradient shrinking Ricci solitons with the stated normalization.
\end{corollary}

The scalar comparison lemma and Theorem~\ref{thm:main} give the lower-bound rigidity theorem. Together with Theorem~\ref{thm:normintro}, this describes both sides of the sharp threshold: a global upper norm bound controls the number of ends, while a global lower norm bound forces cylindrical rigidity in the noncompact case.

\subsection{The complementary scalar upper bound}\label{subsec:scalarupper}
We now assume that $(M^n,g,f)$ is connected and complete, satisfies \eqref{eq:normalization}, and obeys
\begin{equation}\label{eq:scalarceiling}
R\le c:=\frac{n-1}{2}\qquad\text{on }M.
\end{equation}

This scalar upper bound implies noncompactness. Write the Ricci eigenvalues in increasing order as $\lambda_1\le\cdots\le\lambda_n$. We first show that the critical set of $f$ is compact and construct a compact core. We then impose additional conditions at critical points to obtain one-endedness.

\begin{proposition}\label{prop:scalaruppercore}
Let $(M^n,g,f)$ be a connected, complete gradient shrinking Ricci soliton satisfying
\[
\Ric+\nabla^2f=\frac12g,\qquad R+|\nabla f|^2=f,
\]
where $R$ is the scalar curvature. Suppose that
\[
R\le c:=\frac{n-1}{2}\qquad\text{on }M.
\]
Then $M$ is noncompact, the critical set $\Crit(f)=\{x\in M:\nabla f(x)=0\}$ is compact, and $M$ has finite topological type and finitely many ends. For every $t>c$, the value $t$ is regular, and
\begin{equation}\label{eq:uppercollar}
M\setminus\{f<t\}\cong\{f=t\}\times[t,\infty).
\end{equation}
In particular, the number of ends equals the number of connected components of the compact level set $\{f=t\}$.
\end{proposition}
\begin{proof}
If $M$ were compact, integrating the trace equation would give
\[
0=\int_M\Delta f\dd=\int_M\left(\frac n2-R\right)\dd
\ge\frac12\Vol(M)>0,
\]
a contradiction. Thus $M$ is noncompact and $f$ is proper by Lemma~\ref{lem:potential}. At any critical point, $f=R\le c$. Hence $\Crit(f)$ is a closed subset of the compact set $\{f\le c\}$ and is compact.

The normalization also gives the gradient estimate
\begin{equation}\label{eq:uppergradient}
|\nabla f|^2=f-R\ge f-c.
\end{equation}
There are no critical points on $\{f>c\}$. For $t>c$, the level set $\{f=t\}$ is nonempty, compact and smooth, and the flow of $\nabla f/|\nabla f|^2$ proves \eqref{eq:uppercollar}. The flow exists across each finite slab of $f$ because that slab is compact and contains no critical points. The compact core $\{f\le t\}$ is connected, since it is a deformation retract of connected $M$. Its interior is diffeomorphic to $M$ by reparametrizing the exterior collar together with an interior boundary collar, which proves finite topological type. A compact smooth level set has finitely many components, each providing one product end in \eqref{eq:uppercollar}.
\end{proof}

The compact-core construction shows how a scalar curvature bound yields finite topological type and finitely many ends. It also gives the absence of critical points near infinity assumed in the quantitative results of Bertellotti--Buzano \cite{BBends}. To study connectedness of the exterior level set, we next derive quantitative restrictions on the Ricci spectrum at a critical point when $M$ has multiple ends.

\begin{proposition}\label{prop:scalarupperobstruction}
Let $(M^n,g,f)$ be a connected, complete gradient shrinking Ricci soliton satisfying
\[
\Ric+\nabla^2f=\frac12g,\qquad R+|\nabla f|^2=f,
\]
where $R$ is the scalar curvature. Suppose that
\[
R\le c:=\frac{n-1}{2}\qquad\text{on }M.
\]
Write the eigenvalues of $\Ric$ with respect to $g$ in increasing order as $\lambda_1\le\cdots\le\lambda_n$, and let $\Crit(f)=\{x\in M:\nabla f(x)=0\}$. If $M$ has more than one end, there exists $p\in\Crit(f)$ such that, with $w(p)=c-R(p)\ge0$,
\begin{equation}\label{eq:uppercriticalobstruction}
\lambda_2(p)\ge\frac12,\qquad
\lambda_1(p)\le-w(p),\qquad
|\Ric|^2(p)\ge\frac{n-1}{4}+w(p)^2.
\end{equation}
Consequently, if $R<c$ everywhere and $M$ has more than one end, there is a negative Ricci eigenvalue at the critical point $p$ above, with absolute value at least $c-R(p)$.
\end{proposition}
\begin{proof}
By Proposition~\ref{prop:scalaruppercore}, $M$ is noncompact and $f$ is proper. Since $f=R\le c$ on $\Crit(f)$, this closed set lies in the compact set $\{f\le c\}$ and is compact. If $\lambda_2<1/2$ at every critical point, then $\nabla^2f=\tfrac12g-\Ric$ would have at least two strictly positive eigenvalues there, and Lemma~\ref{lem:compactcore} would give one-endedness. By contraposition, multiple ends force $\lambda_2(p)\ge1/2$ at some critical point. Thus
\[
\sum_{i=2}^n\lambda_i(p)\ge\frac{n-1}{2}=c,
\qquad
\lambda_1(p)=R(p)-\sum_{i=2}^n\lambda_i(p)\le R(p)-c=-w(p).
\]
Since $w(p)\ge0$, we have $\lambda_1(p)^2\ge w(p)^2$. Summing the squares and using $\lambda_i(p)\ge1/2$ for $i\ge2$ gives
\[
|\Ric|^2(p)=\lambda_1(p)^2+\sum_{i=2}^n\lambda_i(p)^2
\ge w(p)^2+\frac{n-1}{4}.
\]
The last assertion follows when $w(p)>0$.
\end{proof}

\subsection{Mixed critical-point conditions for one-endedness}\label{subsec:mixedcritical}

To prove one-endedness, it suffices to exclude at least one of the preceding obstructions at each critical point. We obtain two complementary conditions that, when the scalar curvature is below $c$, allow some negative Ricci curvature or a Ricci norm above the cylindrical threshold.

\begin{proposition}\label{prop:mixedoneend}
Let $(M^n,g,f)$ be a connected, complete gradient shrinking Ricci soliton satisfying
\[
\Ric+\nabla^2f=\frac12g,\qquad R+|\nabla f|^2=f,
\]
where $R$ is the scalar curvature. Suppose that
\[
R\le c:=\frac{n-1}{2}\qquad\text{on }M.
\]
Let $\Crit(f)=\{x\in M:\nabla f(x)=0\}$, and let $\lambda_{\min}(\Ric)(p)$ denote the smallest eigenvalue of $\Ric_p$ with respect to $g_p$. Suppose that at every $p\in\Crit(f)$ at least one of the following two strict inequalities holds:
\begin{align}
\lambda_{\min}(\Ric)(p)&>R(p)-c,\label{eq:mixedRiccritical}\\
|\Ric|^2(p)&<\frac{n-1}{4}+\bigl(c-R(p)\bigr)^2.\label{eq:mixednormcritical}
\end{align}
The choice between the two inequalities may depend on $p$. Then $M$ has exactly one end and $H_{n-1}(M;\mathbb Z_2)=0$. In particular, imposing either one of these inequalities on the entire critical set suffices. No Ricci sign or norm hypothesis away from the critical set is required.
\end{proposition}
\begin{proof}
Suppose $\lambda_2(p)\ge1/2$ at a critical point. The eigenvalue calculation in the proof of Proposition~\ref{prop:scalarupperobstruction} then gives both
\[
\lambda_{\min}(p)\le R(p)-c,
\qquad |\Ric|^2(p)\ge\frac{n-1}{4}+\bigl(c-R(p)\bigr)^2.
\]
Both \eqref{eq:mixedRiccritical} and \eqref{eq:mixednormcritical} therefore fail at $p$, contrary to the hypothesis. Thus $\lambda_2<1/2$ at every critical point, so $\nabla^2f=\tfrac12g-\Ric$ has at least two strictly positive eigenvalues there. Proposition~\ref{prop:scalaruppercore} gives noncompactness and properness of $f$. Moreover, $\Crit(f)$ is closed and $f=R\le c$ on it, so it is contained in the compact set $\{f\le c\}$ and is compact. Lemma~\ref{lem:compactcore} therefore gives one-endedness and the mod-two homology vanishing, also in the presence of degenerate critical points.
\end{proof}

Under $R\le c$, positive Ricci curvature on the critical set implies \eqref{eq:mixedRiccritical}. Proposition~\ref{prop:mixedoneend} therefore gives the following corollary.

\begin{corollary}\label{cor:upperpositivecritical}
Let $(M^n,g,f)$ be a connected, complete gradient shrinking Ricci soliton satisfying
\[
\Ric+\nabla^2f=\frac12g,\qquad R+|\nabla f|^2=f,
\]
where $R$ is the scalar curvature. Suppose that
\[
R\le\frac{n-1}{2}\qquad\text{on }M,
\]
and that $\Ric_p>0$ at every $p\in\Crit(f)=\{x\in M:\nabla f(x)=0\}$. Then $M$ has exactly one end and $H_{n-1}(M;\mathbb Z_2)=0$.
\end{corollary}

\begin{remark}\label{rem:upperstrictrecovery}
Under the global bound $R\le c$, the critical-set assumptions $R<c$ and $\Ric\ge0$ give $\lambda_{\min}\ge0>R-c$ on $\Crit(f)$. Thus the mixed criterion recovers the conclusion of Theorem~\ref{thm:criticalintro} under this additional global scalar bound. More generally, \eqref{eq:mixedRiccritical} permits negative eigenvalues, provided their magnitude is smaller than the scalar deficit $c-R$.
\end{remark}

If Ricci curvature is only nonnegative on the critical set, the strict mixed inequality may fail at the scalar threshold. In this case, multiple ends force the cylindrical Ricci spectrum and a rank-one Hessian at some critical point.

\begin{proposition}\label{prop:upperendpoint}
Let $(M^n,g,f)$ be a connected, complete gradient shrinking Ricci soliton satisfying
\[
\Ric+\nabla^2f=\frac12g,\qquad R+|\nabla f|^2=f,
\]
where $R$ is the scalar curvature. Suppose that
\[
R\le c:=\frac{n-1}{2}\qquad\text{on }M,
\]
and that $\Ric_p\ge0$ at every $p\in\Crit(f)=\{x\in M:\nabla f(x)=0\}$. If $M$ has more than one end, there exists a critical point $p$ with
\begin{equation}\label{eq:upperendpointspectrum}
R(p)=c,\qquad
\operatorname{spec}(\Ric_p)=\left\{0,\frac12,\ldots,\frac12\right\},
\end{equation}
where $\operatorname{spec}(\Ric_p)$ denotes the eigenvalues of $\Ric_p$ with respect to $g_p$, counted with multiplicity. At this point $\nabla^2f$ has one eigenvalue $1/2$ and $n-1$ zero eigenvalues. Consequently, if $\operatorname{rank}(\nabla^2f)_p\ne1$ at every $p\in\Crit(f)$, then $M$ has exactly one end. In particular, the additional assumption that $f$ is Morse suffices.
\end{proposition}
\begin{proof}
Choose the critical point in Proposition~\ref{prop:scalarupperobstruction}. Ricci nonnegativity there and $\lambda_i\ge1/2$ for $i\ge2$ give
\[
c\ge R(p)=\lambda_1(p)+\sum_{i=2}^n\lambda_i(p)\ge c.
\]
Equality forces $\lambda_1=0$ and $\lambda_2=\cdots=\lambda_n=1/2$, proving \eqref{eq:upperendpointspectrum}. The Hessian statement follows from $\nabla^2f=\tfrac12g-\Ric$. By contraposition, excluding these rank-one critical points gives one-endedness. In particular, a Morse critical point has a full-rank Hessian and cannot have this exceptional spectrum.
\end{proof}

\begin{remark}\label{rem:upperendpointgap}
Equation~\eqref{eq:upperendpointspectrum} gives a necessary condition at a single critical point for multiple ends under the stated hypotheses. A global rigidity conclusion, such as constant scalar curvature or a cylindrical structure, requires a further argument under $R\le c$ and Ricci nonnegativity only on the critical set when exceptional rank-one critical points are allowed.
\end{remark}

\end{document}